\documentclass{elsarticle}
\usepackage{amsmath}
\usepackage{amssymb}
\usepackage{graphics}
\usepackage{epstopdf}
\usepackage{float}
\newtheorem{theorem}{Theorem}[section]

\begin{document}

\begin{frontmatter}

\title{BIFURCATION ANALYSIS OF WILSON-COWAN MODEL WITH SINGULAR IMPULSES}

\author[a]{Marat ~Akhmet\corref{cor}}
\ead{marat@metu.edu.tr}

\author[a]{Sabahattin ~\c{C}a\u{g}}
\ead{sabahattincag@gmail.com}

\cortext[cor]{Corresponding author}

\address[a]{Department of Mathematics, Middle East Technical University, Ankara, Turkey.}

\begin{abstract}
The paper concerns with Wilson-Cowan neural model with impulses. The main novelty of the study is that besides the traditional singularity of the model, we consider singular impulses. A new technique of analysis of the phenomenon is suggested. This allows to consider the existence of solutions of the model and  bifurcation in ultimate neural behavior is observed through numerical simulations. The bifurcations are reasoned by impulses and singularity in the model and they concern the structure of attractors, which consist of newly introduced sets in the phase space such that medusas and rings.
\end{abstract}


\end{frontmatter}


\section{Introduction}

Wilson and  Cowan \cite{wandc} proposed a model  for describing the dynamics of localized populations of excitatory and inhibitory neurons. This model is a coarse-grained description of the overall activity of a large-scale neuronal network, employing just two differential equations \cite{Kilpatrick}. It is used in the developing of multi-scale mathematical model of cortical electric activity with realistic mesoscopic connectivity \cite{epilepsy}.
On the other hand, sudden changes and the instantaneous perturbations in a neural network at a certain time, which are identified by external elements, are examples of impulsive phenomena which may influence the evolutionary process of the neural network \cite{akca}. In fact, the existence of impulse is often a source of richness for a model. That is to say, the impulsive neural networks will be  an appropriate description of
symptoms of sudden dynamic changes. Therefore, the models considered in this paper have impulsive moments.

The singularly perturbed problems depend on a small positive parameter, which is in front of the derivative, such that the solution varies rapidly in some regions and varies slowly in other regions. They  arise in the various processes and phenomena such as chemical kinetics, mathematical biology, neural networks, fluid dynamics and in a variety models for control theory \cite{segel,hek,owen,terman,fluid, kokotovic,gondal}. In this article,  we will investigate the Wilson-Cowan model with singular impulsive function in which singular perturbation method has been used to analyze the dynamics of neuronal models.

Local bifurcations are ubiquitous in mathematical biology \cite{kuang1} and mathematical neuroscience \cite{hop,diego,dias}, because they provide a framework for understanding behavior of the biological networks modeled as dynamical systems. Moreover, a local bifurcation can affect the global dynamic behavior  of a neuron \cite{hop}. There are many neuronal models  to consider the bifurcation analysis, for instance, the bifurcation for Wilson-Cowan model is discussed in the book of Hoppensteadt and Izhikevich \cite{hop} in which they consider the model of the following type 
\[\dot{x}=-x+S(\rho+cx),\]
where $x\in \mathbb{R}$ is the activity of the neuron,  $\rho\in \mathbb{R}$ is the external input to the neuron,  the feedback parameter  $c\in \mathbb{R}$ characterizes the non-linearity of the system, and $S$ is a sigma shaped function. This system consists only one neuron or one population of neurons. When the bifurcation parameter $\rho$ changes, the saddle-node bifurcation occurs. In our paper, we will discuss two and four of population of  neurons. These systems have impulses of prescribed moments of time. We will observe the local bifurcation in these models. 

The attractors observed in our simulations do not resemble any attractors which have already been observed in the literature. This is why, we need to introduce a new terminology to describe an ultimate behavior of motion in the model.
We call the recently introduced components of constructed attractors as medusas and rings.  
This ``zoological" approach to dynamics is not unique in differential equations. For example, canards 
are cycles of singularly perturbed differential equations \cite{krupa,szmolyan,de}. They were discovered in the van der Pol oscillator  by Benoit et al \cite{benoit}. This phenomenon explains the very fast transition upon variation of a parameter from a small amplitude limit cycle to a relaxation oscillation \cite{krupa}. 
The fast transition is called canard explosion and happens within an exponentially small range of the control parameter. Because this phenomenon is hard to detect it was nicknamed a canard, after the French newspaper slang word for hoax. Furthermore, the shape of these periodic orbits in phase space resemble a duck; hence the name ``canard," the French word
for duck. So the notion of a canard cycle was born and the chase after these creatures began \cite{campbell}. It is important to note that both canards and medusas appear in the singularly perturbed systems.

Bifurcation occurred in this paper cannot be reduced to the existing local bifurcations in the literature, namely, saddle-node, pitchfork, Hopf bifurcations, etc. First of all, we are talking about the change of an attractor set in the four subpopulations of neurons of Wilson-Cowan model with impulses depending on the change of the small parameter. This time the bifurcation parameter is also the parameter of the singularity. Moreover, it is a parameter of the singularity not only in the differential equations of the model, but also in the impulsive part of it. Thus, the cause of bifurcation is not the change of eigenvalues, but it relates to the singular compartment and the impulsive dynamics of the model. This is why, theoretical approvement of the observed bifurcations has not been done in the paper. However, we see that the abrupt changes in the phase portrait through simulations. Additionally, we notice that in the numerical study attractors of the model can be described through the new picture's elements which we call as medusa, medusa without ring and rings, which, in general,  may not be considered invariant for solutions of the model despite that the elements are introduced for the first time. We are confident that they are very generic for differential equations with impulses and they will give a big benefit in the next investigations of discontinuous neural networks.

We will start by defining the membrane time constant since it will be used as the parameter of singularity and bifurcation.

\section{Membrane Time Constant}

\begin{figure}[H]
	\centering
	\includegraphics[scale=0.07]{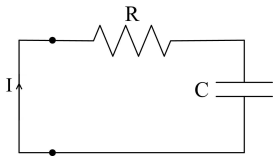}
	\vspace{-10pt}
	\caption{A simple RC circuit.}
	\label{fig:RC}
\end{figure}

The role of the membrane time constant is important in Wilson-Cowan models. In these models the frequency of the oscillation is determined primarily by the membrane time constants \cite{ramesh}. Let us define the membrane time constant $\mu$ for a simple circuit.
Suppose that the membrane is characterized by a single membrane capacitance $C$  in series with a single voltage-independent membrane resistance $R,$ see Fig. \ref{fig:RC}. 
Then, by Ohm's law  the dynamics of the potential $V$ across this circuit in response to a current injection $I$ changes as
\begin{equation*}
RC\frac{dV}{dt}=-V+IR,
\end{equation*}
which has the solution \[V(t)=IR(1-e^{-\frac{t}{RC}})\]
The membrane time constant, here, is defined by the product of the membrane resistance and membrane capacitance $\mu=RC.$ The potential $V(t)$ is governed by exponential decay toward
the steady-state $V=IR$ as $\mu \to 0.$ 
The membrane time constant  is used  to understand how quickly a neuron's voltage level changes after it receives  an input signal.

\section{Singular Model with Singular Impulsive Function}

The dynamics of excitatory and inhibitory neurons are described as follows \cite{wandc}
\begin{equation}\label{model}
\begin{split}
\mu_e\frac{dE}{dt}&=-E+(k_e-r_eE)S_e(c_1E-c_2I+P), \\
\mu_i\frac{dI}{dt}&=-I+(k_i-r_iI)S_i(c_3E-c_4I+Q),
\end{split}
\end{equation} 
where $E(t)$ and $I(t)$ are the proportion of excitatory and inhibitory cells firing per unit time at time $t$, respectively,  $c_1$ and $c_2$ are the connectivity coefficients, which are both positive, represent the average number of excitatory and inhibitory synaptic inputs
per cell, $P(t)$ represents the external input to the excitatory subpopulation, the quantities $c_3,c_4$ and $Q(t)$ are defined  similarly for the inhibitory subpopulation. The nonzero quantities $\mu_e$ and $\mu_i$ represent the membrane time constants while $k_e,k_i,r_e$
and $r_i$ are associated with the refractory terms. Moreover, $S_e(x)$ is the sigmoid function  of the following form 
\begin{equation}
S_e(x)=\frac{1}{1+\exp[-a_e(x-\theta_e)]}-\frac{1}{1+\exp(a_e\theta_e)},
\end{equation}
where $\theta_e$ is the position of the maximum slope of $S_e(x)$ and $\max[\dot{S}_e(x)] = a_e / 4,$ and $S_i$ is defined similarly. 

Since the external inputs influence the neurons activities, $E(t)$ and $I(t)$ can change abruptly. It is natural to consider the previous continuous dynamics in the way that the membrane time constants proceed to be involved in the electrical processes and the impulsive equations have the form
\begin{equation}\label{impulsefunction2}
\begin{split}
&\Delta E|_{t=\theta_i}=\bar{K}(E,I), \\
&\Delta I|_{t=\theta_i}=\bar{J}(E,I),
\end{split}
\end{equation}
where the impulse moments $\theta_i$s are distinct, $\theta_i \in (0,T)$ and the equality $\Delta E|_{t=\theta_i} = E(\theta+)-E(\theta-)$ denotes the jump operator in which $t=\theta$ is the time when the external input influence $E(t)$, $E(\theta-)$ is the pre-impulse value and $E(\theta+)$ is the post-impulse value. Moreover, if one considers the impulsive equations as the limit cases of the differential equations, then at some moments impulsive changes of the activities can depend on the membrane time constants, similar to the ones for the system \eqref{model}. More precisely, we will also study the equations of the form 
\begin{equation}\label{impulsefunction1}
\begin{split}
\mu_e\Delta E|_{t=\eta_j}=K(E,I,\mu_e), \\
\mu_i\Delta I|_{t=\eta_j}=J(E,I,\mu_i),
\end{split}
\end{equation}
where the moments $\eta_j$s are, in general, different from $\theta_i$s. Finally, gathering all the dynamics details formulated above, our single Wilson-Cowan model with impulses has the following form
\begin{equation}\label{modelg}
\begin{split}
&\mu_e\frac{dE}{dt}=-E+(k_e-r_eE)S_e(c_1E-c_2I+P), \\
&\mu_i\frac{dI}{dt}=-I+(k_i-r_iI)S_i(c_3E-c_4I+Q),\\
&\Delta E|_{t=\theta_i}=\bar{K}(E,I), \\
&\Delta I|_{t=\theta_i}=\bar{J}(E,I), \\
&\mu_e\Delta E|_{t=\eta_j}=K(E,I,\mu_e), \\
&\mu_i\Delta I|_{t=\eta_j}=J(E,I,\mu_i),
\end{split}
\end{equation} 
with the initial activity $(E(0),I(0))=(E_0,I_0).$

Define the function $F(E,I)=\begin{pmatrix}
-E+(k_e-r_eE)S_e(c_1E-c_2I+P) \\
-I+(k_i-r_iI)S_i(c_3E-c_4I+Q)
\end{pmatrix}.$

Suppose that $E,I  \in \mathbb{R}$, $t \in [0, T],$ $F(E,I)$ is continuously differentiable on $D$, $K(E,I,\mu_e),J(E,I,\mu_i)$  are continuous on $D\times[0,1]$ and $\bar{K}(E,I),\bar{J}(E,I)$  are continuous on $D$,  $D$ is the domain $D=\{0\leq t\leq T, |E|<d,|I|<d \},$ $\theta_i,i=1,2,\dots,p,$ and $\eta_j,j=1,2,\dots,\bar{p},$ are distinct discontinuity moments in $(0,T).$

Substituting $\mu_e=\mu_i=0$ in \eqref{model} and \eqref{impulsefunction1}, we obtain
$F(E,I)=0$ and 
\begin{equation}\label{impulsefunction1mu0}
\begin{split}
&0=K(E,I,0), \\
&0=J(E,I,0).
\end{split}
\end{equation}
Assume that equations $F(E,I)=0$ and \eqref{impulsefunction1mu0} have the steady states \[(E_1,I_1), (E_2,I_2),...(E_k,I_k),(E_{k+1},I_{k+1}),\dots,(E_l,I_l)\] such that all of them  are real and isolated in $\bar{D}$. They are considered to be states of low level background activities since such activities seem ubiquitous in neural tissue. $E(t)$ and $I(t)$ will be used to refer the activities in the respective subpopulations.

The following condition are required for system \eqref{model}. 
\begin{itemize}
	\item[(C1)] Jacobian matrices of $F(E,I)$ at the points $(E_1,I_1), (E_2,I_2),...,(E_k,I_k)$ are Hurwitz matrices (they have eigenvalues whose real parts are negative).
\end{itemize}
This condition implies that the  states $(E_1,I_1), (E_2,I_2),...,(E_k,I_k)$ are stable steady states of the differential equation  \eqref{model}. Moreover, for the impulsive functions we need the following conditions.
\begin{itemize}
	\item[(C2)] For each $j \in \{1,2,\dots,k\}$ there exists $i \in \{1,2,\dots,k\}$ such that
	\[
	\begin{pmatrix}
	E_j \\ 
	I_j
	\end{pmatrix} 
	+  	
	\begin{pmatrix}
	\bar{K}(E_j,I_j) \\ 
	\bar{J}(E_j,I_j)
	\end{pmatrix}
	=
	\begin{pmatrix}
	E_i \\ 
	I_i
	\end{pmatrix} .
	\]
\end{itemize}
That is, after the each impulse moment $\theta_j$ the activity $(E(t),I(t))$ will be close to another stable steady state $\begin{pmatrix}  E_i \\   I_i  \end{pmatrix}.$ 
\begin{itemize}
	\item[(C3)] \[ \lim_{\substack{(E,I)\to (E_j,I_j) \\ \mu_{e,i} \to 0}}
	\begin{pmatrix}
	\frac{K(E,I,\mu_{e})}{\mu_{e}} \\ 
	\frac{J(E,I,\mu_{i})}{\mu_{i}}
	\end{pmatrix}=
	\begin{pmatrix}
	0 \\ 0
	\end{pmatrix}
	, \quad j=1,2,\dots,k.\]
\end{itemize}
In the denominator of the limit we have small parameters $\mu_e$ and $\mu_i$ which decay to zero. In order to avoid a blow up we need the last condition. In addition, the zero value of the limit gives us the privilege that the activities stay in the domain of attractions of the stable steady states.

Denote $D_j$ as the domain of attraction of stable steady state  $(E_j,I_j),j=1,2,\dots,k,$ such that $D_i\cap D_j=\emptyset$ if $i\neq j$ and $D_j\subset D,j=1,2,\dots,k.$ Also, $z_j(t)$ will be used for denoting the solution of $F(E,I)=$ and \eqref{impulsefunction1mu0} 
such that if the initial value $(E_0,I_0)\in D_j,$ then $z_j(t)=(E_j,I_j)$ for $t\in (0,\theta_1]$ and it alternates to the other stable steady states by condition (C2) for the next intervals $(\theta_i,\theta_{i+1}],i=1,2,\dots,p-1.$
\begin{theorem} \label{thml} 
	Suppose that conditions (C1)-(C3) are true. If the initial value $(E_0,I_0)$ is located in the domain of attraction $D_j$  of the steady state $(E_j,I_j),j=1,2,\dots,k,$ then  the solution $(E(t),I(t))$ of \eqref{modelg} with $(E_0,I_0)$ exists on $[0,T]$ and it is satisfies the limit
	\begin{equation}\label{lim}
	\lim_{\mu_{e,i} \to 0}(E(t),I(t))=z_j(t) \quad \text{for} \quad 0<t\leq T,
	\end{equation}
	where $j=1,2,\dots,k(k-1)^p.$
\end{theorem}
The proof follows from the proof in \cite{cag}.

\begin{figure}[H]
	\centering
	\vspace{-10pt}
	\includegraphics[scale=0.25]{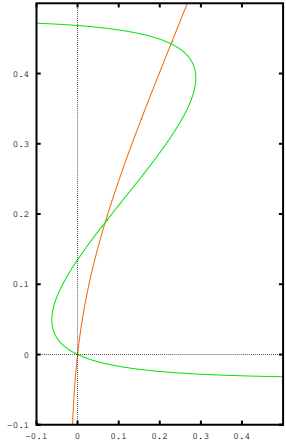}
	\vspace{-10pt}
	\caption{E-I phase plane of \eqref{model0}. The green represents $-E+(0.97-E)S_e(12E-4I)=0$ and the red  represents $-I+(0.98-I)S_i(13E-11I)=0.$}
	\label{fig:wlc}
\end{figure}

\textbf{Example. }  
Now, let us take the external forces $P(t)=Q(t)=0$, $\mu_e=\mu_i=\mu,$ and other coefficients in \eqref{model} as follows: $c_1=12, c_2=4, c_3=13, c_4=11, a_e=1.2, a_i=1,\theta_e=2.8,\theta_i=4, r_e=1, r_i=1, k_e=0.97, k_i=0.98.$ Then, one obtains
\begin{equation}\label{model0}
\begin{split}
\mu\frac{dE}{dt}&=-E+(0.97-E)S_e(12E-4I), \\
\mu\frac{dI}{dt}&=-I+(0.98-I)S_i(13E-11I).
\end{split}
\end{equation}
Taking $\mu=0,$ one has the three equilibria (see Fig. \ref{fig:wlc}), namely
\[
\begin{pmatrix}
E \\ 
I
\end{pmatrix} 
=
\begin{pmatrix}
0 \\ 
0
\end{pmatrix} 
,
\begin{pmatrix}
0.44234 \\ 
0.22751
\end{pmatrix} 
,	
\begin{pmatrix}
0.18816 \\ 
0.067243
\end{pmatrix}. 
\] 

We have $F(E,I)=\begin{pmatrix}
-E+(0.97-E)S_e(12E-4I) \\
-I+(0.98-I)S_i(13E-11I)
\end{pmatrix}.$ Then, the Jacobian matrices of $F(E,I)$ on the steady states are
\[
\begin{pmatrix}
-0.5468 & -0.1511\\ 
0.2250 & -1.1904
\end{pmatrix} 
,
\begin{pmatrix}
-0.9895 & -0.2829\\ 
2.1299 & -31045
\end{pmatrix} 
,\]
\[
\begin{pmatrix}
1.0001 & -0.7469\\ 
0.9879 & -1.9096
\end{pmatrix}, 
\]
respectively. All eigenvalues of the first two matrices are negative, but last one has a positive eigenvalue. Therefore, the first two steady states are stable.

\begin{figure}[H]
	\centering
	\includegraphics[scale=0.7]{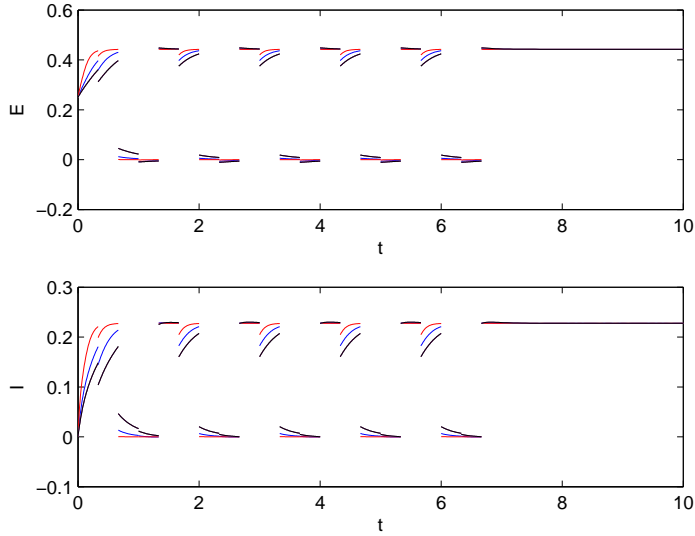}
	\vspace{-10pt}
	\caption{Coordinates of \eqref{model0},\eqref{impulse1},\eqref{impulse2} with initial value $(0.25,0),$ where red, blue and black lines corresponds to value of $\mu=0.1,0.2,0.3,$ respectively.}
	\label{fig:wcex}
\end{figure}

We extend model \eqref{model0} with the following impulse functions
\begin{equation}\label{impulse1}
\begin{split}
&\Delta E|_{t=\theta_i}=-2E+0.44234, \\
&\Delta I|_{t=\theta_i}=-2I+0.22751.
\end{split}
\end{equation}
\begin{equation}\label{impulse2}
\begin{split}
&\mu \Delta E|_{t=\eta_i}=-\mu E^{1/2}(E-0.44234)^2-\sin(\mu^2)I, \\
&\mu \Delta I|_{t=\eta_i}=-\mu I^{1/3}(I-0.22751)^3-\sin(\mu^2)E, 
\end{split}
\end{equation}
where $\theta_i=\frac{2i}{3}, \eta_i=\frac{2i-1}{3}, i=1,2,\dots,20.$ Let us check the conditions of Theorem \ref{thml}. We have shown that the states 
$\begin{pmatrix}
0 \\ 
0
\end{pmatrix} 
,
\begin{pmatrix}
0.44234 \\ 
0.22751
\end{pmatrix}$ 
are stable. Moreover, they satisfy the equations \eqref{impulse2} if $\mu=0.$ Condition (C2) holds since
\[
\begin{pmatrix}
0 \\ 
0
\end{pmatrix} 
+  	
\begin{pmatrix}
0.44234\\ 
0.22751
\end{pmatrix}
=
\begin{pmatrix}
0.44234\\ 
0.22751
\end{pmatrix}
\]
and
\[
\begin{pmatrix}
0.44234\\ 
0.22751
\end{pmatrix} 
+  	
\begin{pmatrix}
-0.88468+0.44234\\ 
-0.45502+0.22751
\end{pmatrix}
=
\begin{pmatrix}
0\\ 
0
\end{pmatrix}.
\]
Lastly, let us check the condition (C3):
\[ \lim_{\substack{(E,I)\to (E_j,I_j) \\ \mu \to 0}}
\begin{pmatrix}
-E^{1/2}(E-0.44234)^2-\frac{1}{\mu}\sin(\mu^2)I\\
-I^{1/3}(I-0.22751)^3-\frac{1}{\mu}\sin(\mu^2)E
\end{pmatrix}=
\begin{pmatrix}
0 \\ 0
\end{pmatrix}
, \quad j=1,2.\]
Clearly, all conditions are satisfied. Therefore, if the initial value $(E_0,I_0)$ is in the domain of attraction of the steady state $(0,0)$ then the activities $(E(t),I(t))$ approaches to the steady states as $\mu \to 0$, that is to say,
\[
\lim_{\mu \to 0}(E(t,\mu),I(t,\mu))=\begin{cases} (0,0) \text{ if } t\in (0,\theta_1]\cup(\theta_2,\theta_3]\cup \dots\\ (0.44234,0.22751) \text{ if } t\in (\theta_1,\theta_2]\cup(\theta_3,\theta_4]\cup \dots \end{cases} ,
\]
and if it is in the domain of attraction of the steady state $(0.44234,0.22751),$ then
\[\lim_{\mu \to 0}(E(t,\mu),I(t,\mu))=\begin{cases} (0.44234,0.22751) \text{ if } t\in (0,\theta_1]\cup(\theta_2,\theta_3]\cup \dots\\ (0,0) \text{ if } t\in (\theta_1,\theta_2]\cup(\theta_3,\theta_4]\cup \dots \end{cases}.\]
To demonstrate the results via simulation, we take  $(E_0,I_0)=(0.25,0)$ which is in the domain of attraction of $(0.44234,0.22751).$ Obviously, the results of the theorem can be seen in Fig. \ref{fig:wcex}.

\section{Bifurcation of New Attractor Composed of Medusa}

In discontinuous dynamics, we  show  that  a new type of attractor consisting medusa, medusa without ring, and rings exist. 
For this purpose,   we study a pair of coupled Wilson-Cowan models. Here, we have four subpopulation in which each system has an excitatory and an inhibitory subpopulation. The first system admits three stable steady states and it has singular impulses. The second one has a limit cycle and it does not have any impulse effects. Also, in the latter system, the membrane time constants are equals to 1.

The first Wilson-Cowan model with impulsive singularity is of the following form:
\begin{equation}\label{3states}
\begin{split}
&\mu_e\frac{dE}{dt}=-E+(0.97-E)S_e(13E-4I), \\
&\mu_i\frac{dI}{dt}=-I+(0.98-I)S_i(22E-2I),\\
&\Delta E|_{t=\theta_i}=6.741E^2-3.58612E+0.45064, \\
&\Delta I|_{t=\theta_i}=6.6087I^2-3.85682I+0.49, \\
&\mu_e \Delta E|_{t=\eta_i}=-\mu_e E^{1/2}(E-0.44234)^2-\sin(\mu_e^2)I, \\
&\mu_i \Delta I|_{t=\eta_i}=-\mu_i I^{1/3}(I-0.22751)^3-\sin(\mu_i^2)E,
\end{split}
\end{equation}
where the sigmoid functions are
\[
\begin{split}
&{S_e}(x)=\frac{1}{1+\exp[-1.5(x-2.5)]}-\frac{1}{1+\exp(3.75)},\\ 
&{S_i}(x)=\frac{1}{1+\exp[-6(x-4.3)]}-\frac{1}{1+\exp(25.8)}.
\end{split}
\]
and impulse moments are $\theta_i=2i+4.95, \eta_i=2i-1+4.95, i=1,2,\dots,50.$ The differential equations in \eqref{3states} have three stable states
\[
\begin{pmatrix}
0 \\ 
0
\end{pmatrix}
,
\begin{pmatrix}
0.20353 \\ 
0.18691
\end{pmatrix} 
,
\begin{pmatrix}
0.45064 \\ 
0.49
\end{pmatrix} 
\]
and two unstable steady states 
\[
\begin{pmatrix}
0.096205 \\ 
0
\end{pmatrix}
,
\begin{pmatrix}
0.37647 \\ 
0.49
\end{pmatrix}. 
\]
The second model, which has a  limit cycle , is of the form 
\begin{align}\label{periodic}
\begin{split}
&\frac{de}{dt}=-e+(0.97-e)\tilde{S_e}(16e-12i+1.25), \\
&\frac{di}{dt}=-i+(0.98-i)\tilde{S_i}(15e-3i),\\
\end{split}
\end{align}
where
\[
\begin{split}
&\tilde{S_e}(x)=\frac{1}{1+\exp[-1.3(x-4)]}-\frac{1}{1+\exp(5.2)},\\ &\tilde{S_i}(x)=\frac{1}{1+\exp[-2(x-3.7)]}-\frac{1}{1+\exp(7.4)}.
\end{split}
\]
We couple system \eqref{3states} and \eqref{periodic} as follows
\begin{align}\label{coupled}
\begin{split}
&\mu_e\frac{dE}{dt}=-E+(0.97-E)S_e(13E-4I), \\
&\mu_i\frac{dI}{dt}=-I+(0.98-I)S_i(22E-2I),\\
&\frac{de}{dt}=-e+(0.97-e)\tilde{S_e}(16e-12i+1.25), \\
&\frac{di}{dt}=-i+(0.98-i)\tilde{S_i}(15e-3i),\\
&\Delta E|_{t=\theta_i}=6.741E^2-3.58612E+0.45064, \\
&\Delta I|_{t=\theta_i}=6.6087I^2-3.85682I+0.49, \\
&\mu_e \Delta E|_{t=\eta_i}=-\mu_e E^{1/2}(E-0.44234)^2-\sin(\mu_e^2)I, \\
&\mu_i \Delta I|_{t=\eta_i}=-\mu_i I^{1/3}(I-0.22751)^3-\sin(\mu_i^2)E.
\end{split}
\end{align}

\begin{figure}[H]
	\centering
	\includegraphics[scale=0.5]{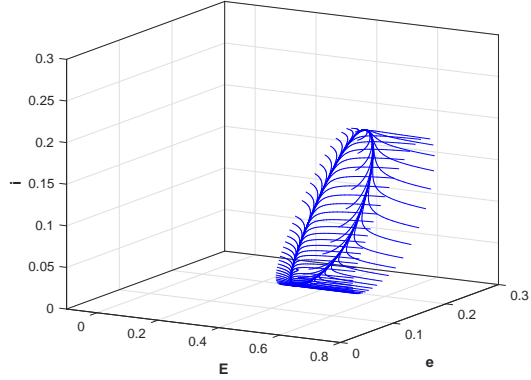}
	\vspace{-10pt}
	\caption{(E,e,i)-coordinates of system \eqref{coupled} for the initial value $(0.4656,0.1101,0.1101,0.04766)$ and the parameter $\mu=0.05.$}
	\label{fig:diadem050}
\end{figure}

It is already known that differential equations in \eqref{3states} has three stable steady states. Suppose that  the membrane time constants in \eqref{coupled} are equal such that $\mu_e=\mu_i=\mu$ and the initial condition is $(0.4656,0.1101,0.1101,0.04766).$ Clearly, in Fig. \ref{fig:diadem050}, one can observe that a  medusa exist for the  value of parameter $\mu=0.05.$ Note that this is a single trajectory and its form looks like a medusa. 

\begin{figure}[H]
	\centering
	\includegraphics[scale=0.5]{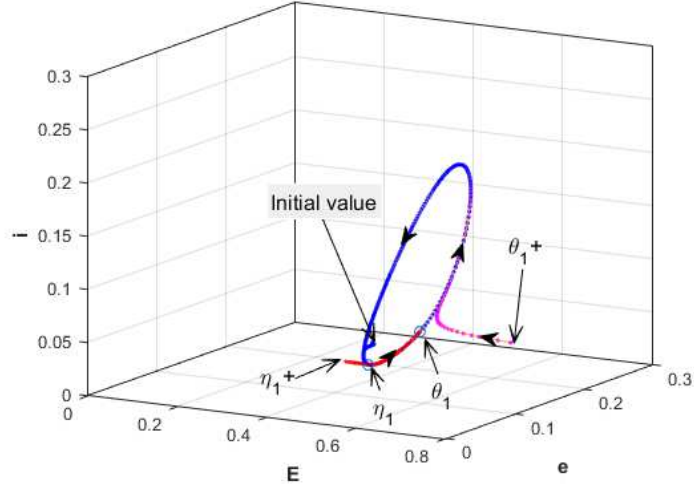}
	\vspace{-10pt}
	\caption{Formation of Fig. \ref{fig:diadem050}. }
	\label{fig:medusa}
\end{figure}

Fig. \ref{fig:diadem050} is formed as follows. The (E,e,i)-coordinate which start at the given initial value approaches to the cycle. It moves around the cycle until the impulse moment $\eta_1.$ When the time reaches  $t=\eta_1,$ because of the impulse function the coordinate jump to $(E(\eta_1+,\mu),e(\eta_1+,\mu),i(\eta_1+,\mu)).$ Again it will approach the cycle and move until the impulse moment $t=\theta_1.$ Then the coordinate jumps to $(E(\theta_1+,\mu),e(\theta_1+,\mu),i(\theta_1+,\mu))$ and it will approach to the  cycle.  The (E,e,i)-coordinate moves in this pattern and finally the medusa in Fig. \ref{fig:diadem050} is observed. The pattern is visualized in Fig. \ref{fig:medusa}.

\begin{figure}[h]
	\centering
	\includegraphics[scale=0.5]{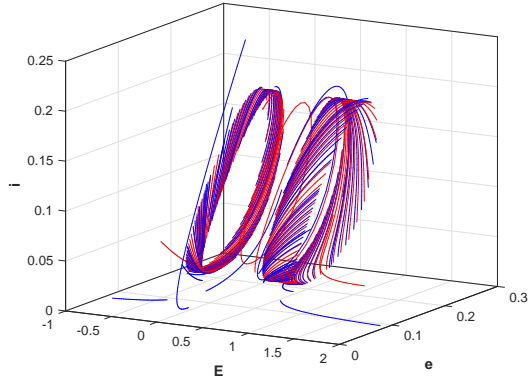}
	\vspace{-10pt}
	\caption{(E,e,i)-coordinates of system \eqref{coupled} for the initial values $(-0.01,0,0.17,0.25),$ $(0.21,0.20,0.20,0.15),$  and for the parameter $\mu=1.$ Blue and red  trajectories are correspondence to the each initial value, respectively. It is seen that two medusas without ring and one cycle are formed. The cycle is between two medusas without ring.}
	\label{fig:diadem10}
\end{figure}

In the following figures, we will see that for different values of the parameter $\mu$ and for the different values of the initial conditions, we will obtain different medusas and rings. First of all,  consider the system \eqref{coupled} with the initial values $(-0.01,0,0.17,0.25),$ $(0.21,0.20,0.20,0.15),$  and with the parameter $\mu=1$ to get Fig. \ref{fig:diadem10}. In this figure, there are two medusas without ring and a cycle. Indeed, they are a single trajectory, which is disconnected in the geometrical sense, but it is connected in the dynamics sense.

Next, we change the parameter to $\mu=0.2$ and use different initial activations $(-0.01,0,0.17,0.25),$ $(0.21,0.20,0.20,0.15),$ $(0.5,0.5,0.3,0.3).$ 

\begin{figure}[H]
	\centering
	\includegraphics[scale=0.5]{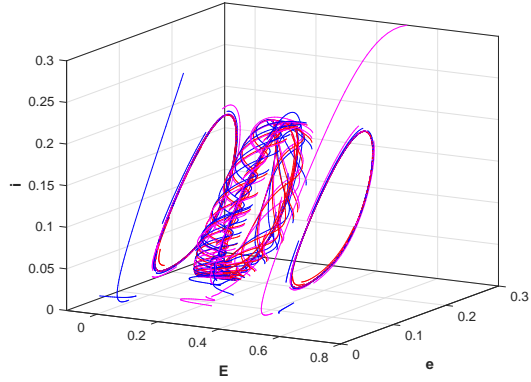}
	\vspace{-10pt}
	\caption{Attractor consists of one medusa and two different rings. Blue, red and magenta trajectories represent solutions in the coordinates (E,e,i) for the given initial values  $(-0.01,0,0.17,0.25),(0.21,0.20,0.20,0.15),(0.5,0.5,0.3,0.3),$ respectively, and  $\mu=0.2.$ 	}
	\label{fig:diadam021}
\end{figure}
In Fig. \ref{fig:diadam021}, one medusa and two different rings are emerged. Geometrically, the attractor is disconnected. However, it is connected in the dynamics sense since it is a single attractor with three parts. There does not exist any limit cycle. The cycles which look like limits cycles are just parts of the whole trajectory.
 
\begin{figure}[H]
 	\centering
 	\includegraphics[scale=0.5]{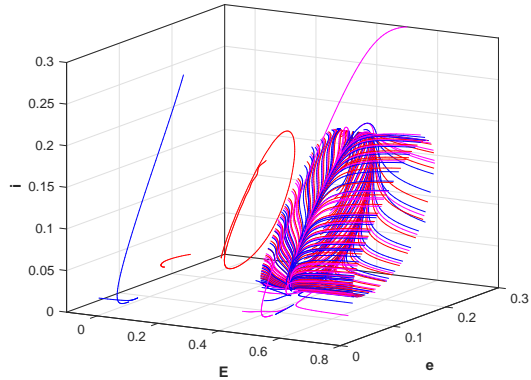}
 	\vspace{-10pt}
 	\caption{Attractor consist only one medusa. Blue, red and magenta trajectories represent solutions in the coordinates (E,e,i) for the given initial values  $(-0.01,0,0.17,0.25),(0.21,0.20,0.20,0.15),(0.5,0.5,0.3,0.3),$ respectively, and  $\mu=0.1.$ The alone red cycle is not an attractor. It is just a part of the trajectory.}
 	\label{fig:diadam011}
 \end{figure}

 Let us consider Fig. \ref{fig:diadam011}. In this figure initial activations are same as in \ref{fig:diadam021}.
 The parameter is fixed and $\mu=0.1.$ Although the initial values are different, the trajectories eventually obtain the shape of the same medusa.  
 There is an alone red trajectory. It is a part of the whole red trajectory. Therefore, it is neither a limit cycle nor a ring since the trajectory never comes to the neighborhood of it.

\begin{figure}[H]
	\centering
	\includegraphics[scale=0.5]{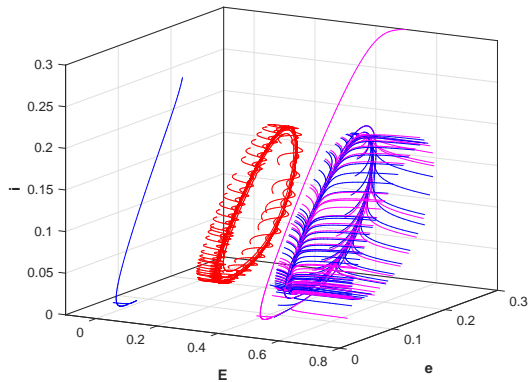}
	\vspace{-10pt}
	\caption{Trajectories of system \eqref{coupled} in coordinates (E,e,i) for different  initial values $(-0.01,0,0.17,0.25),(0.21,0.20,0.20,0.15),
		(0.5,0.5,0.3,0.3)$ and for the fixed  parameter $\mu=0.05.$ Blue, red and magenta trajectories represent solutions for the given initial values, respectively. 
	}
	\label{fig:diadem052}
\end{figure}

Finally, fix the parameter $\mu=0.05.$ In Fig. \ref{fig:diadem052}, any trajectory from the different initial values blue$(-0.01,0,0.17,0.25),$ red$(0.21,0.20,0.20,0.15),$ magenta$(0.5,0.5,0.3,0.3),$ ultimately gets the form of  red or blue medusa. The blue and the magenta trajectories  converges to the same medusa. This is why, we will say that the attractor consists of two disjoint medusas. They are disjoint since there is not a single trajectory which makes two medusas. 

Note that as the parameter decreases the form of the trajectory becomes horizontal through the E-coordinate.

In conclusion, we see that the neuron's dynamics have the following properties: in the case $\mu=1$ two medusas without ring and a cycle are obtained. When $\mu=0.2$ one medusa and two rings emerge and when $\mu=0.1$ one medusa  emerges. Finally, if $\mu=0.05$ two medusas emerge. These results demonstrate that for different values of $\mu$  the qualitative changes in the behavior of trajectories  of \eqref{coupled} occur ultimately. Therefore, we have a bifurcation. It is important to note that
this bifurcation occurs because of the singularity and impulses. This is why, one cannot explain the bifurcations in this paper through the  traditional types of bifurcations, saddle-node, pitchfork, Hopf bifurcation, etc. For example, the change of the numbers of medusas and rings in the local phase portrait depend on the impulsive jumps sizes. Bifurcation, here, also depends on the positions of cycles for unperturbed system.

\section{Conclusion}

It is the first time that only a single small parameter $\mu$ causes not only  to the singularity, but also to the bifurcation. The singularity in this paper is a new kind such that it emerges both from the differential equation part and in the impulsive function. It is also important that the small parameter $\mu$ is a natural parameter which comes from the membrane time constant in Wilson-Cowan neuron model.
We have shown the existence of bifurcation through the simulations. Theoretical proofs are not given since it is difficult to analyze the discontinuous dynamics of the model in which a single parameter  causes both singularity and bifurcation. Therefore, bifurcation is not occurred by the change of eigenvalues, but it relates to the singular compartment and the impulsive dynamics of the model.

New type of attractor, which consists medusa, medusa without ring and rings, is defined. The name comes from the similarity of the form of trajectory and medusa.

\section*{References}

\bibliography{mybibfile}

\begin{thebibliography}{10}
\expandafter\ifx\csname url\endcsname\relax
  \def\url#1{\texttt{#1}}\fi
\expandafter\ifx\csname urlprefix\endcsname\relax\def\urlprefix{URL }\fi
\expandafter\ifx\csname href\endcsname\relax
  \def\href#1#2{#2} \def\path#1{#1}\fi

\bibitem{wandc}
H.~R. Wilson, J.~D. Cowan, Excitatory and inhibitory interactions in localized
  populations of model neurons, Biophysical Journal 12 (1972) 1--24.

\bibitem{Kilpatrick}
Z.~P. Kilpatrick, Wilson-cowan model, in: D.~Jaeger, R.~Jung (Eds.),
  Encyclopedia of Computational Neuroscience, Springer New York, New York,
  2013, pp. 1--5.

\bibitem{epilepsy}
Y.~Wang, M.~Goodfellow, P.~N. Taylor, G.~Baier, Dynamic mechanisms of
  neocortical focal seizure onset, PLoS Comput Biol 10~(8) (2014) 1--18.

\bibitem{akca}
H.~Akça, R.~Alassar, V.~Covachev, Z.~Covacheva, E.~Al-Zahrani, Continuous-time
  additive hopfield-type neural networks with impulses, Journal of Mathematical
  Analysis and Applications 290~(2) (2004) 436 -- 451.

\bibitem{segel}
L.~A. Segel, M.~Slemrod, The quasi-steady state assumption: a case study in
  perturbation, SIAM Review 31 (1989) 446--477.

\bibitem{hek}
G.~Hek, Geometric singular perturbation theory in biological practice, J. Math.
  Biol. 60 (2010) 347--386.

\bibitem{owen}
M.~R. Owen, M.~A. Lewis, How predation can slow, stop, or reverse a prey
  invasion, Bulletin of Mathematical Biology 63 (2001) 655--684.

\bibitem{terman}
D.~Terman, Dynamics of singularly perturbed neuronal networks, in: An
  introduction to mathematical modeling in physiology, cell biology, and
  immunology ({N}ew {O}rleans, {LA}, 2001), Vol.~59 of Proc. Sympos. Appl.
  Math., Amer. Math. Soc., Providence, RI, 2002, pp. 1--32.

\bibitem{fluid}
E.~R. Damiano, R.~D. Rabbitt, A singular perturbation model of fluid dynamics
  in the vestibular semicircular canal and ampulla, Journal of Fluid Mechanics
  307 (1996) 333--372.

\bibitem{kokotovic}
P.~V. Kokotovic, Applications of singular perturbation techniques to control
  problems, SIAM Review 26 (1984) 501--550.

\bibitem{gondal}
I.~Gondal, On the application of singular perturbation techniques to nuclear
  engineering control problems, IEEE Transactions on Nuclear Science 35 (1988)
  1080--1085.

\bibitem{kuang1}
Y.~Kuang, Delay Differential Equations: With Applications in Population
  Dynamics, Mathematics in Science and Engineering, Elsevier Science, 1993.

\bibitem{hop}
F.~C. Hoppensteadt, E.~M. Izhikevich, Weakly Connected Neural Networks, Applied
  Mathematical Sciences, Springer, New York, 1997.

\bibitem{diego}
D.~Fasoli, A.~Cattani, S.~Panzeri, The complexity of dynamics in small neural
  circuits, PLoS Computational Biology 12~(8) (2006) e1004992.

\bibitem{dias}
A.~P.~S. Dias, J.~S. Lamb, Local bifurcation in symmetric coupled cell
  networks: Linear theory, Physica D: Nonlinear Phenomena 223~(1) (2006) 93 --
  108.

\bibitem{krupa}
M.~Krupa, P.~Szmolyan, Relaxation oscillation and canard explosion, Journal of
  Differential Equations 174~(2) (2001) 312--368.

\bibitem{szmolyan}
P.~Szmolyan, M.~Wechselberger, Canards in r3, Journal of Differential Equations
  177~(2) (2001) 419 -- 453.

\bibitem{de}
P.~De~Maesschalck, N.~Popovi{\'c}, T.~J. Kaper, Canards and bifurcation delays
  of spatially homogeneous and inhomogeneous types in reaction-diffusion
  equations, Advances in Differential Equations 14~(9-10) (2009) 943--962.

\bibitem{benoit}
E.~Benoit, J.~L. Callot, F.~Diener, M.~Diener, et~al., Chasse au canard
  (premi{\`e}re partie), Collectanea Mathematica 32~(1) (1981) 37--76.

\bibitem{campbell}
S.~A. Campbell, E.~Stone, T.~Erneux, Delay induced canards in a model of high
  speed machining, Dynamical Systems 24~(3) (2009) 373--392.

\bibitem{ramesh}
R.~Srinivasan, S.~Thorpe, P.~Nunez, Top-down influences on local networks:
  Basic theory with experimental implications, Frontiers in Computational
  Neuroscience 7 (2013) 29.

\bibitem{cag}
M.~Akhmet, S.~\c{C}a\u{g}, A differential equation with singular impulses and
  multi-stable roots, Submitted.

\end{thebibliography}

\end{document}